
\magnification=\magstep1

\input amstex

\documentstyle{amsppt}

\leftheadtext{E. Makai, Jr.}

\rightheadtext{The volume product problem}

\topmatter

\title The recent status of the volume product problem\endtitle

\author E. Makai, Jr.* \endauthor

\address Alfr\'ed R\'enyi Mathematical Institute, 
\newline Hungarian Academy of Sciences,
\newline H-1364 Budapest, Pf. 127, 
\newline HUNGARY
\newline {\rm{http://www.renyi.hu/\~{}makai}}
\endaddress

\email makai\@renyi.hu\endemail

\thanks *Research (partially) supported by Hungarian National Foundation for 
Scientific Research, grant nos. K75016, K81146, and the European
Commission Project TODEQ (MTKD-CT-2005-030042)\endthanks


\subjclass {\it 2010 Mathematics Subject Classification.} Primary:
52A40. Secondary: 52A38, 52A20\endsubjclass

\abstract In this small survey we consider the volume product, and sketch some
of the best upper and lower estimates known up to now, based on our paper
[BMMR]. The author thanks the organizers of the conference in Jurata, March
2010, for their kind invitation, and the excellent atmosphere there. This
paper is based on the talk of the author on that conference.
\endabstract

\endtopmatter\document

\head I. Introduction \endhead

For a finite dimensional Banach space $X$ (we always consider only {\it{real
Banach spaces}}), the {\it{volume product of}} $X$ is the
product of the volumes of the unit balls of the space $X$ and its dual $X'$.
(Of course, for this to make sense, 
we have to identify $X$ and $X'$, which can be done 
via the usual scalar product. However, the 
volume product is independent of the scalar product used.)
This concept has
a great importance in the local theory of Banach spaces, i.e., the
asymptotical (functional analytical) 
investigation of Banach spaces of high finite dimension, cf. e.g., Pisier's
book [Pi]. 
Moreover, it lies
on the cross-road of several disciplines of mathematics, even seemingly
unrelated, which points out its great importance and usefulness.

In more geometric terms,
in ${\Bbb R}^n$ 
we have an $0$-symmetric convex body $K$, and we consider its polar
body $K^*$, and the product of their volumes $V(K)V(K^*)$.

Here a set $K \subset {\Bbb R}^n$ is a {\it{convex body}}, if
it is a compact convex set with non-empty interior. We denote by $V(K)$
its volume.

Let us suppose that $0 \in {\text{int}}\,K$. Then the {\it{polar body $K^*$ of
$K$}} is defined as $K^*:= \{ x \in {\Bbb R}^n \mid \forall k \in K \,\,\,\,
\langle x,k \rangle \le 1 \} $. This is also a convex body, with 
$0 \in {\text{int}}\,K^*$. We have $(K^*)^*=K$. If $K$ is $0$-symmetric, hence
is the unit ball of an $n$-dimensional Banach space $X$, then $K^*$ is the
unit ball of the dual space $X'$. (Of course, for this to make sense, 
we have to identify $X$ and $X'$, via the usual scalar product.)  

We call in the above situation (i.e., when 
$0 \in {\text{int}}\,K$) the quantity
$V(K)V(K^*)$ the {\it{volume product of}} $K$.
This is invariant under non-singular linear transformations. 

The investigation
of this quantity was originated by Blaschke [Bl], who used this concept for the
affine geometry of convex bodies (i.e., investigation of properties of convex
bodies that are invariant under --- possibly volume preserving --- {\it{affine 
transformations}}, i.e., maps of the form $x \mapsto Tx+a$, where $T:{\Bbb
R}^n \to {\Bbb R}^n$ is a
non-singular linear map, and $a \in {\Bbb R}^n$ is a vector),
and Mahler [Ma38] and [Ma39], who used this concept
in the geometry of numbers (i.e., a common part of geometry and number theory,
dealing, among others,
with relation of convex bodies, and also of their polars, to
{\it{lattices}} in ${\Bbb R}^n$, i.e., non-singular linear images of ${\Bbb
Z}^n$ in ${\Bbb R}^n$).

Later it became obvious, that the volume product is a very important quantity 
in the local theory of Banach spaces, it has relations to a number of other
characteristics of these Banach spaces (cf., e.g., [Pi]). 
This fact drew the interest of
functional analysts to the volume product problem, which later resulted in the
solution of the problem of the lower estimate of the volume product for an
$n$-dimensional Banach space in a way satisfactory for applications in the 
local theory of Banach spaces. (The problem of the sharp
upper estimate had already been solved by this time.)

Several other mathematical disciplines also need the volume product. We
mention stochastic geometry, i.e., geometric probability (in several different
ways),
and the geometry of Minkowski spaces (these are just finite dimensional Banach
spaces, but the object is not their analytic properties, but their geometric
properties, which are interesting already for the case $n=2$).
It is interesting enough, that 
the theory of functions of several complex variables is
also connected to the volume product problem. Still one area is discrete
geometry, that is the branch of geometry dealing with density estimates
of systems of convex sets in ${\Bbb R}^d$ 
satisfying various hypotheses. The {\it{density}} of such a system
is the ``percentage of the volume of the whole space covered by the system,
taken with multiplicity (i.e., doubly covered parts count doubly, etc.)''.

So the question is: what is the minimum, and the maximum of the volume
product $V(K)V(K^*)$. 
Seemingly a bit more generally (but in fact equivalently): 
for $x \in {\text{int}}\,K$ we consider 
$V(K)V\left( (K-x)^* \right) $, and investigate this quantity. Here, for the
upper estimate one has to be a bit careful, since for ${\text{dist}}\,(x,\,\,
{\text{bd}}\,K) \to 0$ we have $V(K)V\left( (K-x)^* \right) \to \infty $.
Therefore we
have to take rather $\min _{x \in {\text{int}}\,K} 
V(K)V((K-x)^*)$ (for the $0$-symmetric case this is just $V(K)V(K^*)$). This
quantity is affine invariant.

We give some examples. 
For $K$ the Euclidean unit ball, whose volume is denoted by 
$\kappa _n$, we have 
$$
V(K)V(K^*)=\kappa _n ^2=n^{-n}\left( 2e \pi +o(1) \right) ^n = 
n^{-n}\left( 17.0794... +o(1) \right) ^n\,.
$$
If $K$ is the cube $[-1,1]^n$, or its polar, the {\it{regular cross-polytope}}
${\text{conv}}\, \{ \pm e_i \} $ (with $e_i$ the standard unit vectors), then
$$
V(K)V(K^*)=4^n/n!=n^{-n}\left( 4e + o(1) \right) ^n =
n^{-n}\left( 10.8731... + o(1) \right) ^n\,.
$$  
 
For $x \in {\text{int}}\,K$ one has the simple formula
$$ 
V\left( (K-x)^* \right) = \frac{1}{n} \int _{S^{n-1}} (h_K(u)-\langle x,u
\rangle )^{-n} du\,,
$$

\newpage

where $h_K(u):= \max \{ \langle k,u \rangle  \mid k \in K \} $ 
is the {\it{support function of}} $K$, defined for $u \in S^{n-1}$.
This readily implies that the second
differential of $V\left( (K-x)^* \right) $ with respect to 
$x\in {\text{int}}\,K$ 
is a positive definite quadratic
form, hence $V\left( (K-x)^* \right) $
is a strictly convex function of $x \in $ int$\,K$. 
Hence there is a unique point $x \in {\text{int}}\,K$, such that $V\left(
(K-x)^* \right) $ is minimal: this point $x$
is called the {\it{Santal\'o point of}}
$K$, and is denoted by $s(K)$. If $K$ is $0$-symmetric, we have $s(K)=0$.
    
As a further example, for $K$ a simplex, we have that $s(K)$ is the barycentre
of $K$, and
$$
\cases
V(K)V[ \left( K-s(K) \right) ^* ] 
=(n+1)^{n+1}/(n!)^2=n^{-n}\left( e^2+o(1) \right) ^n = \\
n^{-n}\left( 7.3890... + o(1) \right) ^n\,.
\endcases
$$

In geometry it is unnatural to restrict ourselves to $0$-symmetric bodies,
which {\it{is}} natural in the local theory of Banach spaces. Actually, the
importance of the Santal\'o point shows up only in the asymmetric case, and
also leads to proofs, even in the $0$-symmetric case,
which would be hardly guessed if we would restrict ourselves to the
$0$-symmetric case only.

So, the correct question is: what is the minimum, and the maximum of $V(K)
\times $ 
$V[ \left( K-s(K) \right) ^* ]$. This quantity is invariant under affinities.
 
{\it{In all our theorems $K \subset {\Bbb R}^n$ will be a convex body.}}


\head II. Upper bound \endhead

\proclaim{Theorem 1} ([Bl], [San], [SR], [Pe], [Ba], [MP])
We have 
$$
V(K)V[ \left( K-s(K) \right) ^* ] \le \kappa _n ^2\,,
$$
with equality if and only if $K$ is an ellipsoid.
\endproclaim


\head III. Lower bound \endhead

\definition{Conjecture A} (Mahler-Guggenheimer-Saint Raymond, [Ma39], [G],
[SR])
If $K$ is $0$-symmetric, we have
$$
V(K)V[ \left( K-s(K) \right) ^* ] =V(K)V(K^*) \ge 4^n/n!\,,
$$
with equality exactly for the unit balls of the Hanner-Hansen-Lima spaces.
\enddefinition

The {\it{Hanner-Hansen-Lima spaces}} 
are inductively defined from lower dimensions,
by taking, for some decomposition $n_1+n_2=n$ (where $n_i \ge 1$), 
either the $l^1$, or $l^{\infty
}$-sums of the already defined Hanner-Hansen-Lima spaces
in $n_1$ and $n_2$ dimensions (in
other words, we take in ${\Bbb R}^n={\Bbb R}^{n_1}\oplus {\Bbb R}^{n_2}$
for the unit ball either 
the convex hull, or the sum of the 
unit balls of ${\Bbb R}^{n_1}$ and ${\Bbb R}^{n_2}$). The
basis of induction is $n=1$, when there is a unique $1$-dimensional Banach
space, which is by definition a Hanner-Hansen-Lima space. 
The unit balls of these 
spaces
include, among others, the cube $[-1,1]^n$, the regular

\newpage

cross-polytope 
${\text{conv}}\,
\{ \pm e_i \} $, but also many other bodies. 
A {\it{Hanner-Hansen-Lima body}} is the unit ball of a Hanner-Hansen-Lima space.
For all these bodies $K$ we have $V(K)V(K^*)=4^n/n!$. 

\definition{Conjecture B} (Mahler, [Ma39])
We have
$$
V(K)V[ \left( K-s(K) \right) ^* ] \ge (n+1)^{n+1}/(n!)^2\,,
$$
with equality exactly for the simplex.
\enddefinition


\head III.1. General theorems \endhead

\proclaim{Theorem 2A} ([BMi], [K], [N])
If $K$ is $0$-symmetric, we have
$$
V(K)V(K^*) \ge \kappa _n ^2/2^n=n^{-n}\left( e \pi + o(1) \right) ^n =
n^{-n}\left( 8.5397... + o(1) \right) ^n\,.
$$
\endproclaim

\proclaim{Theorem 2B} ([BMi], [K])
We have
$$
V(K)V[\left( K-s(K) \right) ^*] \ge a n^{-n}( e \pi /2 ) ^n = 
n^{-n}\left( 4.2698... + o(1) \right) ^n\,,
$$
for some constant $a>0$.
\endproclaim

We remark that [BMi] proved only $\kappa _n ^2/A^n$, for a non-explicit
constant $A>0$, in both theorems. 
They have used for this aim the so called {\it{quotient of subspace
theorem}} of V. D. Milman (cf., e.g., [Pi])
that is another very important theorem in the 
local theory of
Banach spaces. This says the following. Any $n$-dimensional Banach space $X$
has a subspace $X_1$, and $X_1$ has a quotient $X_2$, such that
${\text{dim}}\,X_2 \ge cn$, 
and such that the Banach-Mazur distance of $X_2$ and
the Hilbert space of dimension dim$\,X_2$ is at most $C$. Here $c \in (0,1)$,
and $C \in (1, \infty )$ are some constants. More exactly: for any 
$c \in (0,1)$ there exists a $C=C(c)\in (1, \infty )$, such that the previous
statement is true, and conversely, for any  
$C \in (1, \infty )$ there exists a $c=c(C) 
\in (0, 1)$, such that the previous
statement is true. Moreover, [N] proved only $(4^n/n!)( \pi /4)^{3n}$, and for
this aim he used the theory of functions of several complex variables.


\head III.2. Sharp theorems for special bodies \endhead

\head III.2.A. Bodies with high symmetry \endhead 

\proclaim{Theorem 3A}
([SR],  [Me86], [R87])
If $K$ is symmetric with respect to all coordinate hyperplanes (thus,
in particular, is $0$-symmetric), then 
$$
V(K)V(K^*) \ge 4^n/n!\,,
$$
with equality exactly for the unit balls of the Hanner-Hansen-Lima spaces.
\endproclaim

\newpage

Bodies with the property in the hypothesis 
of this theorem, and also the corresponding norms,
are called {\it{unconditional}}.

\proclaim{Theorem 3B} ([BF])
If $K$ has all the symmetries (i.e., congruences) of a regular simplex, then
$$
V(K)V[\left( K-s(K) \right) ^*] \ge (n+1)^{n+1}/(n!)^2\,,
$$
with equality exactly for a regular simplex.
\endproclaim


\head III.2.B. Zonoids \endhead

A {\it{zonoid}} is a limit, in the Hausdorff-metric, of some sequence of
finite sums of segments.

\proclaim{Theorem 4} ([R85], [R86], [GMR]) 
Let $K$ be a convex body, that is also an $0$-symmetric zonoid. Then
$$
V(K)V(K^*) \ge 4^n/n!\,,
$$
with equality if and only if $K$ is a parallelepiped, with centre at \,$0$.
\endproclaim

Of course, this also means that the other Hanner-Hansen-Lima 
bodies are not zonoids.


\head III.2.C. Planar case \endhead

\proclaim{Theorem 5A} ([Ma38], [R86])
Let $n=2$ and let $K$ be $0$-symmetric. Then
$$
V(K)V(K^*) \ge 8\,,
$$
with equality if and only if $K$ is a parallelogram with centre at \,$0$.
\endproclaim

\proclaim{Theorem 5B} ([Ma38], [Me91])
Let $n=2$. Then
$$
V(K)V[ \left( K-s(K) \right) ^*] \ge 27/4\,,
$$
with equality if and only if $K$ is a triangle.
\endproclaim


\head III.2.D. Local minima \endhead

\proclaim{Theorem 6A} ([NPRZ])
Among $0$-symmetric bodies $K$, the volume product $V(K)$ 
$\times V(K^*)$ has a strict
local minimum for $0$-symmetric parallelepipeds.
\endproclaim

\proclaim{Theorem 6B} ([KR])
The volume product $V(K)V[ \left( K-s(K) \right) ^*]$ has a strict
local minimum for simplices.
\endproclaim


\head III.2.E. Polyhedra with small numbers of vertices, or of
$(n-1)$-faces \endhead

\newpage

\proclaim{Theorem 7A} ([LR])
Let $n \le 8$, and let $K$ be an $0$-symmetric convex
polyhedron, with at most $n+1$
opposite (i.e., symmetric w.r.t. $0$)
pairs of vertices, or of \,$(n-1)$-faces. Then
$$
V(K)V(K^*) \ge 4^n/n!\,,
$$
with equality if and only if $K$ is a Hanner-Hansen-Lima body, 
with at most $n+1$
opposite pairs of vertices, or of \,$(n-1)$-faces.
\endproclaim

\proclaim{Theorem 7B} ([MR])
Let $K$ be a convex polyhedron, with at most $n+3$
vertices, or $(n-1)$-faces. Then
$$
V(K)V[ \left( K-s(K) \right)^*] \ge (n+1)^{n+1}/(n!)^2\,,
$$
with equality if and only if $K$ is a simplex.
\endproclaim

Polyhedra (parallelepipeds) in higher dimensional Euclidean
spaces are usually called {\it{polytopes (parallelotopes)}} in geometry.


\head IV. Stability variants \endhead

If one proves an inequality, and also determines the cases of equality, that
is not yet the end of the story. One can further 
ask the following. If a convex body
has a volume product $V(K)V[ \left( K-s(K) \right) ^*]$ (or $V(K)V(K^*)$ for
the $0$-symmetric case)
$\varepsilon $-close to the extremal value, does it
follow, that $K$ is $\left( 1+f( \varepsilon ) \right) $-close, in the
Banach-Mazur distance, to some of the extremal
bodies. Here $f$ is some positive function, 
with $\lim _{ \varepsilon \to 0} f(\varepsilon )=0$.
(Of course, in case of compactness, 
the existence of such a
function $f$ is  evident, but the aim is to give an explicit such function.) 
Typically, at such theorems, $f( \varepsilon )$ is some constant times some
power of $\varepsilon $, sometimes also with some 
logarithmic factor, also on some power. 
At such theorems it is considered as satisfactory, if the sharp order of
magnitude of the function $f$ is determined.

The {\it{Banach-Mazur distance}} $\delta _{BM}(K, L)$
of two $0$-symmetric 
convex bodies $K$ and $L$ is the Banach-Mazur distance of the Banach spaces
with these unit balls. In geometrical terms, this can be described as
follows.
$$
\cases
\delta _{BM} (K,L) = \min \{ \lambda _2 / \lambda _1 \mid \lambda _1, \lambda
_2>0,\,\,
\exists \,\,T
{\text{ non-singular linear map,}} 
\\
{\text{such that }} 
\lambda _1 K \subset TL \subset
\lambda _2 K \} \,.
\endcases
$$
One can extend this naturally to the not $0$-symmetric case as
follows. {\it{For $K$ and $L$ convex bodies,}}
$$
\cases
\delta _{BM} (K,L) := \min \{ \lambda _2 / \lambda _1 \mid \lambda _1, \lambda
_2>0,\,\,
\exists \,\,T
{\text{ non-singular linear map,}}
\\
\exists \,\,a_1,a_2 {\text{ vectors, such that }}
\lambda _1 K +a_1 \subset TL \subset \lambda _2 K +a_2 \} \,.
\endcases
$$
For $K$ and $L$ both $0$-symmetric, this reduces to the previous formula.
(We recall that actually $\log \delta _{BM}(K,L)$ is a metric.)

\newpage

Then, for Theorem 1, 
we have stability, but probably not with the sharp order 
([Bo], [BB], [BMa]). 
Concerning the theorems with lower estimates, we have stability of 
Theorem 4 ([BH]), Theorem 5A ([BH],[BMMR]), 
of Theorem 5B ([BMMR]), of Theorem 6A 
([NPRZ]), and of Theorem 6B ([KR]), at each of these with the sharp order, 
namely, with linear order.


\head V. Functional variants \endhead

Recently there have appeared several functional variants of the volume product
problem. We discuss here only one of these. The objects of investigation are 
the {\it{log-concave functions}}, 
i.e., functions $f:{\Bbb R}^n \to [0, \infty )$, whose
logarithm is concave. To a convex body $K$, with $0 \in {\text{int}}\,K$, we
associate the function $f(x):={\text{exp}}\,(- \| x \| _K ^2 /2)$, 
where 
$\| \cdot \|_K$ means the asymmetric norm (i.e., $\| -x \| _K 
\ne \| x \| _K$) associated to the ``unit ball'' $K$ ({\it{i.e.,}} 
$\| x \| _K := \min \{ \lambda \ge 0 \mid x \in \lambda K \} $).
Then $V(K)={\text{const}}_n 
\cdot \int _{{\Bbb R}^n} f(x)\,dx$, so 
here the right hand side is the proper generalization of $V(K)$. 

Moreover, polarity
of $K$ and $K^*$ goes over to the following. The corresponding functions $f$
and $f^*$ have their negative logarithms, with values in 
$(- \infty , \infty ]$,
which are the Legendre transforms of
each other. The {\it{Legendre transform}} of a function $\varphi :{\Bbb R}^n
\to [ -\infty , \infty ]$ is the function ${\Cal L} \varphi :{\Bbb R}^n
\to [ -\infty , \infty ]$, defined by 
$$
({\Cal L}
\varphi )(y):= \sup \{ \langle x,y \rangle - \varphi (x) \mid x \in {\Bbb R}^n
\} \,.
$$ 
Thus, the object of investigation is 
$$
\int _{{\Bbb R}^n} f(x)\,dx \, \cdot \int _{{\Bbb R}^n} f^*(x)\,dx\,,
$$ 
where one supposes
$$
\int _{{\Bbb R}^n} f(x)\,dx \in (0, \infty )\,.
$$
Cf. the nice exposition in [A-AKM].

Unfortunately, translations of convex bodies have no (good) generalizations
to log-concave functions. 
Thus, in place of a translation $K \mapsto K-x$, where $x \in {\text{int}}\,K$,
one considers an arbitrary translate of the function $f$ (i.e., $x \mapsto 
f(x-x_0)$), and proves the sharp upper bound for a suitable translate of the
original function $f$. Here, 
for even functions $f$, one may choose $x_0=0$ (as for
$0$-symmetric bodies one may choose $x=s(K)=0$), cf. [A-AKM].
Of course, the problem of translations 
does not concern the question of the
lower bound (as it is a minimum problem), 
but, in case of the upper bound, only the $0$-symmetric case of the
volume product problem generalizes this way to even log-concave functions.
For the upper bound, in the even case, the sharp upper bound for the
functional variant 
(cf. [Ba], [A-AKM], [FM07])
immediately implies
the sharp upper bound for the
$0$-symmetric case of the volume product problem.
Namely, the extremal
even functions (up to constant factors) 
are ones associated to $0$-symmetric convex bodies, 

\newpage

more exactly, 
to $0$-symmetric
ellipsoids. For
the lower bound, the best, known lower bound for the
functional variant (cf., [FM08a], Theorem 7)
implies the best, known lower bound for the volume
product problem, apart from the actual value of the base of the exponential
(when the second lower bound is written in the form $n^{-n}\left( \right.
$const$\left. +o(1) \right) ^n $). 
In the case of {\it{unconditional 
functions}} (i.e., $f(x_1,...,x_n)=f(|x_1,...,|x_n|)$), with $f>0$,
the sharp lower bound is known (cf., [FM08a], Theorem 6).

Still we note that the conjecture in ${\Bbb R}^n$ 
about the lower bound for the
functional variant, for the even, or the general case (cf., [FM08b]), 
would imply Conjecture A or Conjecture 
B, about the lower bound for the volume product, in the
$0$-symmetric, or the general case, for ${\Bbb R}^n$, or ${\Bbb R}^{n-1}$,
respectively (cf., [FM08b]). 
However, the
conjecture about the lower bound for the
functional variant, for the even, or the general case, {\it{for all}} $n$, is 
equivalent to Conjecture A, or Conjecture B, 
{\it{for all}} $n$, respectively (cf. [FM08b], Propositions 1, 2).


\enddocument